# BROWNIAN EXCURSIONS AND PARISIAN BARRIER OPTIONS: A NOTE

Michael Schröder

*Lehrstuhl Mathematik III*
*Seminargebäude A5, Universität Mannheim, D–68131 Mannheim*

**Abstract.** This paper addresses the Paris barrier options of [**CJY**] and their valuation using the Laplace transform approach. The notion of the Paris barrier options is extended such that their valuation is possible at any point during their lifespan, and the pertinent Laplace transforms of [**CJY**] are modified when necessary.

**1. Introduction:** An interesting application of the theory of Brownian excursion theory has been observed in [**CJY**] in connection with so–called barrier options. These are widely traded options which come as puts or calls that are activated or deactivated as soon as their underlying hits a prespecified barrier level. As explained in [**CJY**], this abruptness entails a number of undesirable consequences and problems both on a theoretical and practical level. The idea of the Paris barrier options of [**CJY**] is to cushion this abruptness. This is by introducing a systematic delay for the consequences of hitting the barrier to become effective. And the insight of [**CJY**] is that well adapted means for this are afforded by excursion theory in general and the Brownian meander and Azéma's martingale in particular.

Still, Paris barrier options by now seem to be traded by the leading institutions only and not in large quantities. And the author is grateful to Pliska (UIC, Chicago) for pointing out to him a principal difficulty with the concepts of [**CJY**]. Trading Paris barrier options will typically lead to situations where the underlying has been staying above or below the respective barrier for a certain period of time but where this excursion has not yet lasted long enough for action to be taken. And there seem to be no results in [**CJY**] how to value Paris barrier options in such situations. In fact, such situations do not seem to be covered by the very constructions in [**CJY**]. One aim of this note so is to provide the necessary conceptual extensions to meet this criticism.

For valuing Paris options, Laplace transforms of the Paris option densities have been computed in [**CJY**]. Pursuing Pliska's suggestions we moreover found it necessary to modify one of these two Laplace transforms. So I am very grateful to Pliska for his remark. And I also wish to thank Yor for the interest he showed in this work and for making available to me corrections to [**CJY**] collected by Jeanblanc–Picqué which have also been incorporated in the exposition of this paper.

**2. Black–Scholes framework:** The analysis of this paper is in the Black–Scholes framework and uses the risk–neutral approach to valuating contingent claims. In this set–up there two securities. First there is a riskless security, a bond, that has the continuously compounding positive interest rate $r$. Then there is a risky security whose price process $S$ is modelled as follows. Start with a complete probability space equipped with the standard filtration of a standard Brownian motion on it that has the time set $[0, \infty)$. On this filtered space one has the *risk neutral measure* $Q$, a probability measure equivalent to the given one, and a standard $Q$–Brownian motion $B$ such that $S$ is the strong solution of the following stochastic differential equation:

$$dS_t = (r - \delta) \cdot S_t \cdot dt + \sigma \cdot S_t \cdot dB_t, \qquad t \in [0, \infty).$$

Herein the positive constant $\sigma$ is the volatility of $S$. The constant $\delta$ depends on the security modelled. For instance, it is the dividend rate if $S$ is a stock.

**3. Paris options:** This section takes up the Paris barrier options of [**CJY**, §2]. Their new idea for cushioning the impact of the underlying hitting the barrier is as follows. They require the underlying $S$ to spend a minimum time $D > 0$ above or below their prespecified barrier $L \geq 0$ before the option is knocked in or knocked out.

The *Paris down–and–in call* to be considered in the sequel is the following European–style contingent claim on $S$ written at time $t_0$ and with time to maturity $T$. Its payoff at $T$ is that of a call on $S$ with exercise price $K$:

$$(S_T - K)^+ = \max\{S_T - K, 0\},$$

if $S_t$ is less than $L$ during a connected subperiod of length at least $D$ of the monitoring period $[t_0, T]$. Equivalently, there is a point in time $a$ such that the interval $I = (a, a+D)$ is contained in the lifespan $[t_0, T]$ of the option and $S_t < L$ for any $t$ in $I$. Otherwise the call expires worthless.

Now let $t$ be any time in $[t_0, T]$ at which the Paris down–and–in call has not yet been knocked in. The time–$t$ value $C_{d,i}$ of this call is then a function of $K$, $\tau = T - t$, $S_0$, $r$, $\delta$, $L$ and $D$, and mostly any of these arguments is suppressed from the notation. The valuation problem depends on an excursion time $H_{L,t}$ defined in the sequel such that the Paris down–and–in call is knocked in before its maturity $T$ if and only if $H_{L,t} \leq T$. Granting this, using the arbitrage pricing principle, the time–$t$ price $C_{d,i}$ of the Paris down–and–in call is then given by the conditional $Q$–expectation:

$$C_{d,i} = e^{-r\tau} E_t\Big[\phi(S_T) \cdot \mathbf{1}_{\{H_{L,t} \leq T\}}\Big],$$

with $\phi(S_T)$ equal to the call payoff $(S_T - K)^+$, and with $\mathbf{1}_{\{H_{L,t} \leq T\}}$ the characteristic function of the event $H_{L,t} \leq T$.

In defining $H_{L,t}$ at time $t$ two constellations are to be distinguished. Namely, either today's security price $S_t$ is equal to or above the barrier $L$, or it is below $L$. In the first case with $S_t \geq L$, the Paris down–and–in call is knocked in before its maturity $T$ iff $S$ remains smaller than $L$ for all points in time of an interval of length at least $D$ contained in $[t, T]$. Thus let $H_{L,t}$ in this case denote the first time $s$ greater than or equal to $t + D$ such that $S_u < L$ for any time $u$ in the interval $(s - D, s)$.

In the second case where $S_t < L$, the security price $S$ has already been staying below $L$ for some connected period time during the lifetime of the option. The down–and–in call is then knocked in before its time to maturity $T$ if the following happens. The security price $S$ continues to stay below $L$ also during the period of time from today until the point in time $t+d$ with $d$ smaller than $D$, i.e., one has $S_u < L$ for all $u$ in $(t, t+d)$. Thus define $H_{L,t} = t+d$ in this case. If the level $L$ is hit by $S$ before time $t+d$, however, the clock for the minimum lenght $D$ is restarted. So define $H_{L,t}$ in this case to be the first time $s$ greater than or equal to $t+D$ such that $S_u < L$ for any $u$ in the subinterval $(s-D, s)$ of the positive real line.

The discussion of Paris barrier options would not be complete without having referred to [**CJY**, §2] for details on the whole family of Paris barrier options. Calls and puts with either of the following barrier types: down–and–out, down–and–in, up–and–out, up–and–in. Valuing these eight types of options is reduced in a standard way to the case of the above Paris down–and–in call. Indeed, put–call parities are given in [**CJY**, §5], and [**CJY**, §4.2] reduces the valuation of the out–calls to that of the in–calls, mutatis mutandis. This paper thus concentrates on valuing the Paris down–and–in call, in the sequel also referred to as *the Paris option* for simplicity.

**4. The basic valuation identity:** Valuing the Paris option in the Black–Scholes context is made precise by the following *basic valuation identity* for its time–$t$ price $C_{d,i}$

$$C_{d,i} = e^{-\left(r+\frac{\varpi^2}{2}\right)\cdot\tau} \int_{\beta(S_t)}^{\infty} e^{\varpi x} \cdot \left(S_t \cdot e^{\sigma x} - K\right) \cdot h_b(\tau, x)\, dx,$$

with the *Paris option density* $h_b$ given by:

$$h_b(u, y) = \int_{-\infty}^{\infty} E^*\left[\mathbf{1}_{\{H_b^* < u\}} \frac{1}{\sqrt{2\pi \cdot (u - H_b^*)}} \cdot e^{-\frac{(x-y)^2}{2\cdot(u-H_b^*)}}\right] d\mu^*(x),$$

for any reals $u > 0$ and $y$, where $\varpi = \sigma^{-1} \cdot (r-\delta-\sigma^2/2)$, with $\beta(S_t) = \sigma^{-1}\log(K/S_t)$, where $H_b^* = H_{L,t} - t$, and with $\mu^*$ the measure for Brownian motion at time $H_b^*$. The expectation is taken with respect to a Girsanov transformed risk neutral measure $Q^*$ constructed in §6 when proving the basic valuation identity. Moreover notice that by construction, $h_b(u,\ )$ is zero for $u$ less than $D$, respectively $d$, depending on today's situation.

**5. Results on the Laplace transforms of the densities:** This section discussed he Laplace transforms with respect to time of the Paris option density $h_b$ of §4. For any fixed real number $y$, recall they are defined by:

$$\mathscr{L}\left(h_b(\ , y)\right)(z) = \int_0^{\infty} e^{-zt} h_b(t, y)\, dt,$$

for any complex number $z$ with sufficiently big real part. Moreover define the function $\Psi$ by:

$$\Psi(z) = \int_0^{\infty} x \cdot e^{-\frac{x^2}{2}+zx}\, dx,$$

for any complex number $z$, and choose on the complex plane with the non–positive real axis deleted the square root defined using the principal branch of the logarithm. The computations of [**CJY**, §§5, 8] then essentially give the following two results.

**Proposition A:** *For any real number $y$, the Laplace transform of $h_b(\,,y)$ is a holomorphic function on the right complex half–plane.*

**Proposition B:** *Suppose $b = \sigma^{-1}\log(L/S_t)$ is non–positive. For any fixed real number $y$, the Laplace transform of $h_b(\,,y)$ is given by:*

$$\mathscr{L}\big(h_b(\,,y)\big)(z) = \frac{e^{\frac{b}{\sqrt{D}}\sqrt{2Dz}}}{\sqrt{D}\sqrt{2Dz}\,\Psi(\sqrt{2Dz})} \int_0^\infty x \cdot e^{-\frac{x^2}{2D} - |b-x-y|\sqrt{2z}}\, dx\,,$$

*for any complex number $z$ with positive real part.*

If $b = \sigma^{-1}\log(L/S_t)$ is positive, i.e., when today's security price $S_t$ is below the barrier $L$, we have modified the Paris option of [**CJY**]. Today's situation is then such that the price of the security has already been staying below the barrier $L$ for a certain connected period of time during the lifetime of the Paris options. For the Paris option to be knocked in the price thus has to stay below the barrier $L$ only for another connected period time of a length $d < D$ from today, i.e., $S_u < L$ for all $u$ in $[t, t+d]$. In this way, Paris options now can be valued at any point of their monitoring periods. The desirability of this modification has been pointed out to me by Pliska. I regard the following modification of the valuation results of [**CJY**] as a direct result of his suggestions.

**Proposition C:** *Suppose $b = \sigma^{-1}\log(L/S_t)$ is positive. For any real number $y$, the Laplace transform of $h_b(\,,y)$ is given on the right half–plane as the following four–term sum of Laplace transforms:*

$$\mathscr{L}\big(h_b(\,,y)\big) = \mathrm{Erfc}\left(\frac{b}{\sqrt{2d}}\right) \cdot \mathscr{L}\big(h_{b,1}(\,,y)\big) + \frac{1}{\sqrt{2\pi d}} \cdot \mathscr{L}\big(h_{b,2}(\,,y)\big)$$

$$+ \mathrm{Erfc}\left(\frac{b}{\sqrt{2d}}\right)\mathrm{Erf}\left(\frac{b}{\sqrt{2d}}\right) \cdot \mathscr{L}\big(h_{b,3}(\,,y)\big) + \mathrm{Erf}\left(\frac{b}{\sqrt{2d}}\right) \cdot \mathscr{L}\big(h_{b,4}(\,,y)\big)$$

Herein $h_{b,k}$ are the functions on the positive real line times the real line defined by:

$$h_{b,3}(u,y) = \mathbf{1}_{(d,\infty)}(u) \cdot \frac{1}{u-d+D}\left\{\frac{\sqrt{u-d}}{\sqrt{2\pi}} \cdot e^{-\frac{(y-b)^2}{2(u-d)}}\right.$$

$$\left. - \frac{1}{2}\frac{(y-b)\sqrt{D}}{\sqrt{u-d+D}} \cdot e^{-\frac{(y-b)^2}{2(u-d+D)}} \mathrm{Erfc}\left(\frac{(y-b)\sqrt{D}}{\sqrt{2(u-d)(u-d+D)}}\right)\right\},$$

$$h_{b,4}(u,y) = \mathbf{1}_{(d,\infty)}(u)\frac{1}{\sqrt{2\pi u}} \cdot f_{b,u}(y) \qquad \text{where} \qquad f_{b,u}(y) = e^{-\frac{y^2}{2u}} - e^{-\frac{(y-2b)^2}{2u}}\,,$$

for any positive real number $u$ and any real number $y$, while the remaining two functions are defined using their Laplace transforms as follows:

$$\mathscr{L}\big((h_{b,1}(\,,y))\big)(z) = \int_0^d \frac{e^{-zw}}{D\sqrt{2z}\,\Psi(\sqrt{2Dz})} \left(\int_0^\infty x \cdot e^{-\frac{x^2}{2D}-|b-x-y|\sqrt{2z}}\,dx\right)\mu_b(dw),$$

$$\mathscr{L}\big((h_{b,2}(\,,y))\big)(z) = \int_0^d \frac{e^{-zw}}{\sqrt{2z}\,\Psi(\sqrt{2Dz})} \left(\int_{\mathbf{R}} e^{-|x-y|\sqrt{2z}}\,f_{b,d}(x)dx\right)\mu_b(dw),$$

for any positive real number $y$ and any complex number $z$ with positive real part. Herein $\mu_b$ denotes the law of the first passage time to the level $b$, given on the positive real line by

$$\mu_b(dw) = \frac{b}{\sqrt{2\pi}}w^{-3/2}e^{-\frac{b^2}{2w}}.$$

**6.   Proof of the basic valuation identity:**  For proving the basic valuation identity, consider the restarted at time $t$ Brownian motion $B^*(u) = B(t+u)-B(t)$ and the normalized Brownian motion with drift $W^*$ given by

$$W_u^* = \frac{1}{\sigma}\log\Big(\frac{S_{t+u}}{S_t}\Big) = \varpi u + B_u^*, \qquad u \in [t,\infty).$$

By construction we then have $S_{t+u} < L$ iff $W^*(u) < b$ and $H_{L,t} \leq T$ iff $H_b^* \leq \tau = T-t$. Using Strong Markov, thus transcribe the price $C_{d,i}$ of §3 of the Paris option as follows:

$$C_{d,i} = e^{-r\tau}E^Q\Big[\big(S_t \cdot e^{\sigma W^*(\tau)}-K\big)^+ \cdot \mathbf{1}_{\{H_b^* \leq \tau\}}\Big].$$

Now proceed in analogy to [**CJY**, §§4,5]. Using the particular case [**KS**, p.196f] of the Cameron–Martin–Girsanov theorem, change measure from $Q$ to $Q^*$ such that $W^*$ becomes a $Q^*$–Brownian motion. Iterating the $Q^*$–expectation thus obtained for $C_{d,i}$, we have

$$C_{d,i} = e^{-\left(r+\frac{\varpi^2}{2}\right)\cdot\tau}E^*\Big[\mathbf{1}_{\{H_b^* \leq \tau\}}E^*\Big[e^{\varpi W_\tau^*}\big(S_t \cdot e^{\sigma W_\tau^*}-K\big)^+\Big|\mathscr{F}_{H_b^*}\Big]\Big].$$

Apply the strong Markov property of Brownian motion using $H_b^*$ as $\mathscr{F}$–stopping time. The resolvent so obtained further simplifies since the random variables $H_b^*$ and $W^*(H_b^*)$ are independent. With $\mu^*$ the measure for $W^*(H_b^*)$, the basic valuation identity follows.

**7.   An intermediate identity for the proofs of the Laplace transforms:**  With the concepts introduced in §4, we have the following intermediate identity:

**Lemma:**  *For any real number $y$, one has:*

$$\mathscr{L}\big(h_b(\,,y)\big)(z) = \int_{\mathbf{R}} E^*\Big[e^{-zH_b^*}\Big]\cdot\frac{e^{-|x-y|\sqrt{2z}}}{\sqrt{2z}}\mu^*(dx),$$

*for any complex number $z$ in the right complex half–plane.*

For proving the Lemma interchange the Laplace transform with the two exponential integrals defining the Paris option density $h_b(u, y)$ in §4 to get:

$$\mathscr{L}(h_b(\,,y))(z) = \int_{\mathbf{R}} E^*\left[\int_{H_b^*}^{\infty} e^{-zu} \frac{1}{\sqrt{2\pi(u - H_b^*)}} \cdot e^{-\frac{(x-y)^2}{2(u-H_b^*)}}\, du\right] \mu^*(dx).$$

In this last integral successively change variables $w = u - H_b^*$ and separate the expectation from the Laplace transform to get:

$$\int_{\mathbf{R}} E^*\left[e^{-zH_b^*}\right] \cdot \mathscr{L}\left(\frac{1}{\sqrt{2\pi u}} \cdot e^{-\frac{(x-y)^2}{2u}}\right)(z)\, \mu^*(dx).$$

Using e.g. [**D**, Beispiel 8, p.50f] for the Laplace transform then completes the proof.

**8. Computing the Laplace transform in the case $b$ equal to zero:** The computation of the Laplace transforms of §5 for $h_b$ is by reduction to the case $b = 0$ where we have to show the *key relation*:

$$E^*\left[e^{-\frac{z^2}{2}H_0^*}\right] \cdot \Psi(z\sqrt{D}) = \Psi(0) = 1,$$

for any complex number $z$ with positive real part. This section reviews the key ideas of its proof in [**CJY**, §8]. The computations are based on the Azéma martingale and properties of the Brownian meander. The *Azéma martingale* is the martingale on $[0, \infty)$ for the slow Brownian filtration $\mathbf{F}^+$ given for any $t \geq 0$ by:

$$\mu_t = (\operatorname{sgn} W_t^*)\sqrt{t - g_t}.$$

Herein $g_t$ is defined as the supremum over all real numbers $s \leq t$ such that $W^*(s) = 0$. The process given by:

$$m_t(u) = \frac{1}{\sqrt{t - g_t}}\bigl|W^*\bigl(g_t + u(t - g_t)\bigr)\bigr|, \qquad u \in [0, 1]$$

is a *Brownian meander* and is independent of the $\sigma$–subfield $\mathscr{F}^+(g_t)$ of the slow Brownian filtration. Its law is independent of $t$ and given by $x \exp(-x^2/2) \cdot \mathbf{1}_{(0,\infty)}\, dx$. The tautology $W^*(t) = m_t(1) \cdot \mu_t$ thus implies the following identity:

$$E^*\left[e^{zW_t^* - \frac{z^2}{2}t}\,\Big|\,\mathscr{F}^+(g_t)\right] = e^{-\frac{z^2}{2}t} \cdot \Psi(z \cdot \mu_t),$$

for any positive real number $z$. Fixing any such $z$, the crucial point is to exhibit the product $\psi(-z \cdot \mu(H_0^*)) \exp(-z^2 H_0^*/2)$ as a martingale. This is achieved by a boundedness argument showing the argument of the above conditional expectation to be uniformly integrable for any time $t$ up to $H_0^*$. Using another form of the optional stopping theorem it follows that the $Q^*$–expectation of the martingale stopped at $H_0^*$ is equal to the $Q^*$–expectation of it at time zero, whence to $\Psi(0)$ which equals one. With the random variables $H_0$ and $W^*(H_0^*)$ independent, the key relation results. This is an identity between functions holomorphic in particular on the right half–plane. Using the identity theorem, it thus remains valid for any complex number $z$ with positive real part.

**9. Computing the Laplace transform for $b$ non–positive:** Using §7 Lemma the proof of §5 Proposition B reduces to compute the expectation of $\exp(-zH_b^*)$, for any complex $z$ with $\operatorname{Re}(z) > 0$, and the density $\mu^*$. Following [**CJY**, §8.3.3], this is by reduction to the results of §8 using restarting arguments. Indeed, decompose $H_b^*$ as follows:

$$H_b^* = T_b + H_0^{**}.$$

Herein $T_b$ is the first passage time of $W^*$ to the level $b$, and $H_0^{**}$ is defined as follows. It is the smallest non–negative point in time $s$ at which the restarted–at–time–$T_b$ Brownian motion $W^{**}(u) = W^*(T_b+u) - b$ is zero for the first time after having been less than zero for a connected period of time of length at least $D$.

To compute the expectation required, notice that $T_b$ and $H_0^{**}$ are independent random variables. The conditional expectation at time $T_b$ of $\exp(-(w^2/2)H_b^*)$ thus is the product of $\exp(-(w^2/2)T_b)$ and the expectation of $\exp(-(w^2/2)H_0^{**})$. With the law of $H_0^{**}$ equivalent to that of $H_0^*$, this last expectation is given by the key relation of §8. So it is deterministic in particular. Take expectations to time zero of this product. The expectation of $\exp(-(w^2/2)T_b)$ that remains to be computed is the Laplace transform of the law of $T_b$ at $w^2/2$. Using [**D**, Beispiel 8, p.50f] for instance, it is standardly seen to be equal to $\exp(-|b|(2 \cdot (w^2/2))^{1/2})$ and thus to $\exp(+bw)$. With $w = (2z)^{1/2}$ one gets:

$$E^*\left[e^{-z H_b^*}\right] = \frac{1}{\Psi(\sqrt{2Dz})} \; E^*\left[e^{-\sqrt{2z}\, T_b}\right] = \frac{1}{\Psi(\sqrt{2Dz})} \; e^{b\sqrt{2z}}.$$

To determine $\mu^*$ notice the tautology $W^*(H_b^*) = W^{**}(H_0^{**}) + b$. Thus $W^*(H_b^*) \leq x$ iff $W^{**}(H_0^{**}) \leq x - b$. From §8 the law of $W^{**}(H_0^{**})$ is minus that of $m_1$, whence

$$Q^*\big(W^*(H_b^*) \in dx\big) = \mathbf{1}_{(-\infty,b]}(x) \cdot (b-x) \cdot e^{-\frac{(x-b)^2}{2D}} \cdot \frac{dx}{D}.$$

Substitute into §7 Lemma and change variables $w = b - x$ to obtain §5 Proposition B for any complex number $z$ with sufficiently big positive real part. Extend to the right half–plane using the consequence of the key relation of §8 that $\Psi((2Dz)^{1/2})$ has no zeroes there. This completes the proof of §5 Proposition B and the $b$–non–positive part of §5 Proposition A.

**10. Computing the Laplace transform for $b$ positive:** The proof of §5 Proposition C modifies the argument of [**CJY**, §8.3.3]. Using §7 Lemma, the first task is to compute the density $\mu^*$ and the $Q^*$–expectation of $\exp(-zH_b^*)$ for any complex number $z$ with positive real part. The two cases where the level $b$ is hit before or after time $d$ are to be distinguished. Thus decompose the underlying probablity space $\Omega$ into the set $A$ on which the first passage time $T_b$ of $W^*$ to the level $b$ is less than or equal to $d$

$$A = \{\omega \in \Omega \,|\, T_b(\omega) \leq d\}$$

and its complement $\Omega \setminus A$ on which $T_b$ is bigger than $d$. This induces the decomposition:

$$E^*\left[e^{-zH_b^*}\right] = E^*\left[\mathbf{1}_A \cdot e^{-zH_b^*}\right] + E^*\left[\mathbf{1}_{\Omega\setminus A} \cdot e^{-zH_b^*}\right]$$

of the $Q^*$–expectation of the random variable $\exp(-zH_b^*)$. By construction $H_b^* = d$ on the complement of $A$. So the second above summand is $\exp(-zd)Q^*(T_b > d)$. On the set $A$ the level $b$ is reached before the critical time $d$ and the clock for the excursion is reset to zero at the first passage time $T_b$ to the level $b$. Accordingly one has the decomposition:

$$H_b^* = T_b + H_0^{**},$$

with $H_0^{**}$ defined as follows. It is the smallest $s \geq 0$ at which the restarted–at–time–$T_b$ Brownian motion $W^{**}(u) = W^*(T_b+u)-b$ is zero for the first time after having been less than zero for a connected period of time of lenght at least $D$. With this random variable $H_0^{**}$ independent from $T_b$, the conditional expectation of the random variable $\mathbf{1}_A \cdot \exp(-zH_b^*)$ at time $T_b$ thus equals $\exp(-zT_b)$ times the expectation of $\exp(-zH_0^{**})$. The key relation of §8 now applies to $H_0^{**}$ and identifies this last expectation as the reciprocal of $\Psi((2zD)^{1/2})$. The expectation of $\mathbf{1}_A \cdot \exp(-zT_b)$ on the other hand, is given by integrating $\exp(-zw)$ from zero to $d$ against the density $\mu_b$ of $T_b$. Summarizing, it so follows for any complex number $z$ with positive real part:

$$E^*\left[e^{-zH_b^*}\right] = Q^*(T_b > d)e^{-zd} + \frac{\int_0^d e^{-zw}\,\mu_b(dw)}{\Psi(\sqrt{2Dz}\,)}\,.$$

For determining $\mu^*$ decompose the random variable $W^*(H_b^*)$ with respect to $A$:

$$W^*(H_b^*) = \left(W^{**}(H_0^{**})+b\right) \cdot \mathbf{1}_A + W^*(d) \cdot \mathbf{1}_{\Omega\setminus A}\,.$$

Using the independence of $T_b$ and $H_0^{**}$, the law $\mu_A^*(dx) = Q^*\big((W^{**}(H_0^{**})+b) \in dx;\, T_b \leq d\big)$ of the first summand is obtained as the convolution of the laws of $T_b$ and $W^{**}(H_0^{**})$:

$$\mu_A^*(dx) = Q^*(T_b \leq d) \cdot \mathbf{1}_{(-\infty, b]}(x) \cdot (b-x) \cdot e^{-\frac{(x-b)^2}{2D}} \cdot \frac{dx}{D}\,.$$

For the law $\mu_{\Omega\setminus A}^*(dx) = Q^*\big(W^*(d) \in dx;\, T_b > d\big)$ of the second summand notice $T_b > d$ iff $W^*(t) < b$ for all $t \leq d$. Thus the first passage time condition $T_b > d$ transscribes into the condition that the running maximum $\max\{W^*(t)|t \leq d\}$ is smaller than $b$. Using [**H**, p.9],

$$\mu_{\Omega\setminus A}^*(dx) = e^{-\frac{x^2}{2d}} \cdot \frac{dx}{\sqrt{2\pi d}} - e^{-\frac{(x-b)^2}{2D}}\,dx \cdot \frac{dx}{\sqrt{2\pi d}}\,.$$

Applying §5 Lemma, the Laplace transform $\mathscr{L}\big(h_b(\,,y)\big)(z)$ of $h_b(\,,y)$ at $z$ is equal to:

$$\mathscr{L}\big(h_b(\,,y)\big)(z) = E^*\left[e^{-zH_b^*}\right] \int_{\mathbf{R}} \frac{e^{-|x-y|\sqrt{2z}}}{\sqrt{2z}}\,(\mu_A^* + \mu_{\Omega\setminus A}^*)(dx)$$

using that the expectation factor is deterministic and independent of the variable $x$. So it suffices to compute the two integrals that make up the second factor of this product. For the first of these, on substituting for $\mu_A^*$ and changing variables $w = b-x$, one gets:

$$\int_{\mathbf{R}} \frac{e^{-|x-y|\sqrt{2z}}}{\sqrt{2z}}\,\mu_A^*(dx) = \frac{Q^*(T_b \leq d)}{D} \int_0^\infty x\,\frac{e^{-|b-x-y|\sqrt{2z}}}{\sqrt{2z}}\,e^{-\frac{x^2}{2D}}\,dx\,.$$

For the second of these one analogously obtains:
$$\int_{\mathbf{R}} \frac{e^{-|x-y|\sqrt{2z}}}{\sqrt{2z}} \mu^*_{\Omega\setminus A}(dx) = \frac{1}{\sqrt{2\pi d}} \int_{\mathbf{R}} \frac{e^{-|x-y|\sqrt{2z}}}{\sqrt{2z}} \left( e^{-\frac{x^2}{2d}} - e^{-\frac{(x-2b)^2}{2d}} \right) dx.$$

Multiplying out, $\mathscr{L}(h_b(\,,y))(z)$ is the sum of four terms of which each is seen to be a Laplace transform itself. Laplace inversion of this identity so is term by term. To identify the remaining functions $h_{b,3}$ and $h_{b,4}$, recall for instance using [**D**, Beispiel 8, p.50f]:
$$\mathscr{L}^{-1}\left( \frac{1}{\sqrt{2z}} \cdot e^{-\alpha\sqrt{2z}} \right)(u) = \frac{1}{\sqrt{2\pi u}} \cdot e^{-\frac{\alpha^2}{2u}}$$

valid for any positive real numbers $\alpha$ and $u$. The Laplace inverse of the improper integral factor of the above $\mu^*_A$ integral is, on applying Fubini's theorem, equal to:
$$\int_0^\infty x e^{-\frac{x^2}{2D}} \cdot \mathscr{L}^{-1}\left( e^{-zd} \frac{e^{-|b-x-y|\sqrt{2z}}}{\sqrt{2z}} \right)(u)\, dx.$$

Successively apply the shifting theorem for the Laplace transform to take care of the factor $\exp(-zd)$ and the above Laplace inversion formula with $\alpha$ equal to $|b-x-y|$ to arrive at:
$$\frac{1}{D} \cdot \mathbf{1}_{(d,\infty)}(u) \frac{1}{\sqrt{2\pi(u-d)}} \int_0^\infty x \cdot e^{-\left(\frac{x^2}{2D} + \frac{(b-x-y)^2}{2(u-d)}\right)} dx.$$

This integral is seen to be the value of the function $h_{b,3}$ at $u$ and $y$. Analogously, the Laplace inverse at any $u > 0$ of the integral factor of the above $\mu^*_{\Omega\setminus A}$ integral is equal to:
$$\frac{1}{\sqrt{2\pi d}} \mathbf{1}_{(d,\infty)}(u) \frac{1}{\sqrt{2\pi(u-d)}} \int_{\mathbf{R}} e^{-\frac{(x-y)^2}{2(u-d)}} \left( e^{-\frac{x^2}{2d}} - e^{-\frac{(x-2b)^2}{2d}} \right) dx.$$

Again suppressing the details of the computation, it is seen to be equal to the value of the function $h_{b,4}$ at $u$ and $y$. This completes the proof of §5 Proposition C.

# Brownian excursions and Parisian barrier options: a note

**Michael Schröder**

(Mannheim)

January 2002